\documentclass[10pt]{amsart}\usepackage{amssymb} \usepackage{amscd}

\usepackage{amssymb}
\usepackage{amscd}
\usepackage{epsfig}

\title[Motivic Iterated Integrals]{Algebraic Cycles and Motivic generic Iterated
Integrals}
\author{Hidekazu Furusho}
\author{Amir Jafari}

\email{furusho@ias.edu, furusho@math.nagoya-u.ac.jp}
\email{amir@ias.edu, amir@math.duke.edu}
\address{Institute For
Advanced Study, School of Mathematics, Princeton, NJ 08540 USA}

\newcommand{\Q}{\mathbb Q}
\newcommand{\Z}{\mathbb Z}
\newcommand{\C}{\mathbb C}
\newcommand{\T}{\mathcal T}

\newcommand{\h}{\mbox{H}}
\newcommand{\N}{\mathcal N}

\newcommand{\I}{\mathbb I}
\newcommand{\To}{\longrightarrow}
\newtheorem{thm}{Theorem}[section]
\newtheorem{lem}[thm]{Lemma}
\newtheorem{cor}[thm]{Corollary}
\newtheorem{prop}[thm]{Proposition}

{\theoremstyle{definition} \newtheorem{rem}[thm]{Remark}
\newtheorem{eg}[thm]{Example}} {\theoremstyle{definition}
\newtheorem{defn}[thm]{Definition}}

{\theoremstyle{remark} }
\begin{document}

\bibliographystyle{amsalpha+}
\maketitle

\begin{abstract}
 Following \cite{GGL}, we will give a combinatorial framework for
 motivic study of iterated integrals on the affine line. We will
 show that under a certain genericity condition these
 combinatorial objects yield to elements in the motivic Hopf
 algebra constructed in \cite{BK}. It will be shown that the Hodge
 realization of these elements coincides with the Hodge structure
 induced from the fundamental torsor of path of punctured affine
 line.
\end{abstract}
\tableofcontents


\section{Introduction.}

Let $F$ be a field. Using cubical algebraic cycles, in \cite{BK}
Bloch and K\v ri\v z constructed a graded Hopf algebra:
$$\chi_F=\bigoplus_{r\ge 0}\chi_F(r).$$

This Hopf algebra should be the Hopf algebra of framed mixed Tate
motives over $F$. In other words the category of mixed Tate motives
over $F$ (MTM$(F)$ for short) should be equivalent to the category
of graded finite dimensional $\chi_F$ comodules. We take this
description as the definition of MTM$(F)$.

 Our main result is construction of ceratin motivic iterated integrals on the
 affine line. We want to stress that our construction works for {\it{any}} field
$F$. When $F$ is a number field \cite{DG} gives a different
construction of motivic fundamental group (torsor) inside the
category of MTM obtained as the heart of a $t$-structure on the
triangulated category of Voevodsky \cite{V}, as explained in
\cite{Lev1}. This yields to a motivic definition of all iterated
integrals on the affine line.

For elements $s_1,\dots,s_n\in S$, we want to construct a motivic
analogue of the iterated integral:
$$\int_a^b \frac{dt}{t-s_1}\circ\dots\circ\frac{dt}{t-s_n}$$
which will be an element of $\chi_F(n)$ and we denote it by ${\Bbb
I}(a;s_1,\dots,s_n;b)$. The special case
$\I(0;\underbrace{1,0,\dots,0}_n;z)$ which (up to a sign)
represents $\mbox{Li}_n(z):=\sum_{k=1}^{\infty} z^k/k^n$, was
constructed in \cite{BK}. Unfortunately we have to assume a
genericity assumption on the sequence $a,s_1,\dots,s_n,b$ in order
to have an admissible cycle. This assumption is as follows:
\begin{defn}
The sequence $a,s_1,\dots,s_n,b$ is generic if either $a=b$ or the
non-zero terms do not repeat.
\end{defn}
{\bf Remark.} In fact the method will work for certain sequences
which are not generic, examples of such sequences are
$$(0;s_1,\dots,s_n;a)$$
where at most two of $s_i$'s is equal to $1$ and the rest are
zero. This will imply the existence of motivic double
polylogarithm $Li_{n,m}(a)$. In \cite{GGL} it is assumed that all
the $s_i$'s are distinct. With allowing zero to repeat we recover
the expression given in \cite{BK} for the polylogarithm, and
moreover we get expressions for multiple polylogarithm
$Li_{n_1,\dots,n_m}(z_1,\dots,z_m)$ for the values $z_i$ such that
$z_i\dots z_j\ne 1$ for $1\le i\le j\le n$.

We will prove that these elements satisfy the usual properties of
the iterated integrals, namely:
\begin{thm}\label{properties}
Under the above genericity assumption, the elements ${\Bbb
I}(a;s_1,\dots,s_n;b)$ satisfy the following properties:
\begin{enumerate}
\item Triviality: \begin{align*}
{\Bbb I}&(a;b)=1, \\
{\Bbb I}&(a;s_1,\dots,s_n;a)=0.
\end{align*}
\item Shuffle relation: $${\Bbb I}(a;s_1,\dots,s_n;b)\cdot {\Bbb
I}(a;s_{n+1},\dots,s_{n+m};b)=\sum_{\sigma\in Sh(n,m)}{\Bbb
I}(a;s_{\sigma(1)},\dots,s_{\sigma(n+m)};b).$$
\item Path composition: $$\I(a;s_1,\dots,s_n;b)=\sum_{k=0}^n
\I(a;s_1,\dots,s_k,c)\cdot\I(c;s_{k+1},\dots,s_n;b).$$

\item Antipode relation:
$$\I(a;s_1,\dots,s_n;b)=(-1)^n\I(b;s_n,\dots,s_1;a).$$
\end{enumerate}
\end{thm}
Furthermore we show the following formula for the coproduct (ref.
\cite{G1}).
\begin{thm}\label{coproduct}
For a generic sequence, the coproduct of $\I(a;s_1,\dots,s_n;b)$
is given by
$$\sum {\mathbb I}(a;s_{i_1},s_{i_2},\dots,s_{i_k};b)\otimes
\prod _{j=0}^k {\mathbb
I}(s_{i_j};s_{i_j+1},\dots,s_{i_{j+1}-1};s_{i_{j+1}})$$ where the
sum is over all indices $0=i_0<i_1<\dots<i_k<i_{k+1}=n+1$ and
$s_0:=a$ and $s_{n+1}:=b$ and $k=0,1,\dots$.
\end{thm}

The content of the paper is as follows. Section 2 contains a
review of definition of cubical cycle complex and the Hopf algebra
$\chi_F$ of \cite{BK}. Definition of bar complex is reviewed, here
we use a different sign convention from \cite{BK}. Section 3 gives
a combinatorial differential graded algebra (DGA) built out of
rooted decorated trees. The main idea is due to Goncharov
\cite{G3}. A morphism from the generic part of this DGA to the
cubical cycle DGA is given. This extends the definition given in
\cite{GGL}. In section 4 for a given tuple
$(a_0;a_1,\dots,a_n;a_{n+1})$ in a finite subset $S$ of ${\Bbb
A}^1(F)$, we give our proposed definition for motivic iterated
integral $\I({\bf a})$. It is built out of a specific cycle
$\rho({\bf a})$ defined using 3-valent trees . The relation
between 3-valent trees and iterated integrals were first studied
by Goncharov in \cite{G3} and \cite{G1}). The crucial property is
a formula for the differential of this cycle. Theorems
\ref{properties} and \ref{coproduct} are proven in section 5,
using the definitions of the previous section. Sections 6 and 7
involve the important calculation of Hodge realization for the
motivic iterated integrals when the ground field is embedded
inside $\C$. This calculation justifies the name motivic iterated
integral for the elements $\I({\bf a})$ constructed in section 4.
In section 8 and 9 we give some complementary remarks for the
non-generic case.
\\
This work is an extension of \cite{GGL}. We have been informed
that Goncharov, Gangl and Levin are writing a more complete paper
which includes some of the details we cover in the present note.

{\bf Acknowlegments.} The authors would like to thank Pierre
Deligne for explaining patiently some of his ideas on mixed
motives. It is apparent from the introduction that this paper owes
a great deal to the works of Goncharov. The second author wants to
thank him for teaching him and sharing his ideas on matters
discussed in this article and much more, during his Ph.D studies
at Brown University. The authors also want to thank IAS for
providing an ideal environment for research. We specially want to
thank Herbert Gangl for answering our questions regarding the
admissibility conditions on cycles in \cite{BK}. We also thank
Kasra Rafi who helped us to construct the differential on trees
while we were unaware of the definitions of \cite{G3} and
\cite{GGL}.
\section{On the Construction of MTM$(F)$ According to Bloch and K\v ri\v z.}



A differential graded algebra (DGA) with Adams grading is a
bi-graded $\Q$-vector space, ${\mathcal A}=\oplus {\mathcal
A}^n(r)$ where $n\in \Z$ and $r\ge 0$, such that ${\mathcal
A}^n(r)=0$ for $n>2r$ and ${\mathcal A}(0)=\Q$, together with a
product:
$${\mathcal A}^n(r)\otimes {\mathcal A}^m(s)\To {\mathcal
A}^{n+m}(r+s)$$ that makes ${\mathcal A}$ into a graded commutative
algebra with identity (the signs are contributed from the
differential grading and not the Adams grading), and a differential
$$d:{\mathcal A}^n(r)\To {\mathcal A}^{n+1}(r)$$
that satisfies the Leibniz rule.

Following \cite{BK}, we now introduce the cubical cycle complex. Let
$F$ be a field and denote $\Box_F:={\mathbb P}_F^1\backslash\{1\}$.
Then the permutation group on $n$ letters $\Sigma_n$ acts on
$\Box_F^n$ and also we have an action of $(\Z/2)^n$ given by
$$\epsilon\cdot(x_1,\dots,x_n)=(x_1^{\epsilon_1},\dots,x_n^{\epsilon_n})$$
where $\epsilon\in \{1,-1\}^n$.
 Therefor we have an action of $G_n:=(\Z/2)^n\rtimes \Sigma_n$ on
$\Box_F^n$. Let $\mbox{Alt}_n\in\Q[G_n]$ be the element
$$\frac{1}{|G_n|}\sum_{g\in G_n}sgn(g)g$$
where for $g=(\epsilon,\sigma)\in G_n$, sign is defined by $(\prod_i
\epsilon_i)sgn(\sigma)$. We also define a face of $\Box_F^n$ as a
subset defined by setting certain coordinates equal zero or
infinity. We are now prepared to define the DGA with an Adams
grading $\mathcal N={\mathcal N}_F$:
$${\N}^n(r):=\mbox{Alt}_{2r-n} Z(\Box_F^{2r-n},r).$$
Notice that although $r\ge 0$, $n$ can be negative. Here
$Z(\Box_F^n,r)$ denotes the $\Q$-span of the admissible codimension
$r$ subvarieties (i.e. closed and integral sub schemes that
intersect all the faces of codimension $\ge 1$ properly, i.e. in
codimension $r$ ) of $\Box_F^n$. The product structure is given by
$$Z_1\cdot Z_2:=\mbox{Alt}(Z_1\times Z_2).$$
The differential is given by
$$d:\N^n(r)\To \N^{n+1}(r)$$
$$\sum_{i=1}^{2r-n}(-1)^{i-1}(\partial_0^i-\partial_{\infty}^i)$$
where for $c=0,\infty$, $\partial^i_c$ is obtained by the pull-back
of the cycles under the inclusions $\Box_F^{2r-n-1}\hookrightarrow
\Box_F^{2r-n}$ given by letting the $i^{\mbox{th}}$ coordinate equal
to $c$. Using the results of \cite{Bl1},\cite{Bl2} and \cite{Lev2}
it is shown in \cite{BK} that:
$$\mbox{H}^n(\N(r))\cong \mbox{CH}^r(F,2r-n)\otimes \Q\cong K_{2r-n}(F)_{\Q}^{(r)}.$$
Here $\mbox{CH}^r(F,n)$ denotes the Bloch's higher Chow group of
$\mbox{Spec}(F)$ and the superscript $r$ for the K-group denotes the
graded quotient with respect to the Adams filtration. Therefore
Bloch and K\v ri\v z argue that under the assumption that $\N$ is
1-connected in the sense of Sullivan (this is a stronger version of
Beilinson-Soul\'{e} vanishing conjecture) its natural to define the
category of mixed Tate motives over $F$ as graded comodules on the
graded Hopf algebra $\h^0(B(\N))$ which is obtained by the bar
construction.

We now recall the bar construction. For our future need we develop
the theory in a more general setting. Let $A$ be a DGA with Adams
grading and let $M$ be a right DG module over $A$ with a compatible
Adams grading. Let $A^+=\oplus_{r>0}A(r)$. Define
$$B(M,A)=M\otimes\left(\bigotimes^{\bullet} A^+[1]\right).$$
The elements of $B(M,A)$ are denoted by a bar notation as
$m[a_1|\dots|a_r]$ which has a degree equal to
$\deg(m)+\deg(a_1)+\dots\deg(a_r)-r$. The graded vector space
$B(M,A)$ is a differential graded with total Adams grading. The
differential $d$ is given as the sum of two differentials, the
external and internal differentials:
\begin{align*}
d_{ext}(m[a_1|\dots |a_r])=&dm[a_1|\dots|a_r] \\
&+\sum_{i=1}^r Jm[Ja_1|\dots|Ja_{i-1}|da_i|a_{i+1}|\dots|a_r],
\\
d_{int}(m[a_1|\dots|a_r])=&Jm\cdot a_1[a_2|\dots|a_r]
\\
&+\sum_{i=1}^{r-1} Jm\cdot [Ja_1|\dots|Ja_{i-1}|Ja_{i}\cdot
a_{i+1}|\dots|a_r].
\end{align*}
 the operation $J$ is given by $J(a)=(-1)^{\deg(a)-1}a$ on
the homogeneous elements. Note that $\deg(a)-1$ is the degree if $a$
in the shifted complex $A[1]$. If $M$ is the trivial DGA, we have a
coproduct (where the empty tensor is 1 by convention):
$$\Delta([a_1|\dots|a_r])=\sum_{s=0}^r [a_1|\dots|a_s]\otimes
[a_{s+1}|\dots|a_r].$$

Up to this point all the constructions work for a DGA which is not
necessarily commutative. The product is defined only when both $M$
and $A$ are commutative DGA's by the shuffle:
$$m[a_1|\dots|a_r]\cdot m'[a_{r+1}|\dots|a_{r+s}]:=\sum
sgn_{{\bf{a}}}(\sigma)m\cdot
m'[a_{\sigma(1)}|\dots|a_{\sigma(r+s)}]$$ where $\sigma$ runs over
the $(r,s)$-shuffles (an $(r,s)$-shuffle is a permutation $\sigma$
on $1,\dots,r+s$ such that $\sigma^{-1}(1)<\dots<\sigma^{-1}(r)$ and
$\sigma^{-1}(r+1)<\dots<\sigma^{-1}(r+s)$) and the sign is obtained
by giving $a_i$'s weights $\deg(a_i)-1$. These data make
$B(A):=B(k,A)$ into a DGA with an Adams grading. It is also a Hopf
algebra.

Furthermore let $\gamma_A$ be the Lie co-algebra of the
indecomposables of $\h^0B(A)$. Bloch and K\v ri\v z show that if $A$
is cohomologically connected then the DGA given by
$\wedge^{\bullet}\gamma_{A}[-1]$ is the 1-minimal model of $A$ in
the sense of Sullivan. The following definition makes sense for any
field $F$, but the desired Ext groups of mixed Tate motives are
obtained if we assume ${\mathcal N}$ is cohomologically connected
(equivalent to Beilinson-Soul\'{e} vanishing conjecture) and is
quasi-isomorphic to its 1-minimal model, i.e.
$\wedge^{\bullet}\gamma_{\mathcal N}[-1]$.
\\
\begin{defn}\label{motive}(\cite{BK}) The category of mixed Tate motives over $F$ is the
category of finite dimensional graded comodules over
$$\chi_F:=\h^0B(\N_F).$$ More explicitly a mixed Tate motive over $F$ is a
graded finite dimensional $\Q$-vector space $M$ with a linear map
$\nu: M\To M\otimes \chi_F$ such that it respects the grading and
$$(id\otimes \Delta)(\nu(a))=(\nu\otimes id)(\nu(a))$$
$$(id\otimes\epsilon)(\nu(a))=a$$
where $\Delta$ and $\epsilon$ are the co-product and the co-unit of
$\chi_F$. \end{defn}

\section{Trees and Cycles.}
The study of the differential structure (and other algebraic
structures such as coproduct, etc.) on graphs is now an active area
in mathematics, with relations to topology, moduli spaces and
physics.

The graph complex was introduced in the seminal work of Kontsevich
in \cite{K}. The differential given here first appeared in
\cite{G3}. This differs from the differential of Kontsevich
because of the existence of external edges. This section is a
complement to \cite{GGL}.

Let $S$ be a set. We define the DGA ${\mathcal T}_S$ (with an
Adams grading) of rooted $S$-decorated planar trees.

\begin{defn}{\bf (Rooted $S$-decorated planar tree).} A rooted $S$-decorated planar tree is a
connected graph with no loops, such that each edge has exactly two
vertex and no vertex has degree 2. Each vertex of degree one is
decorated by an element of the given set $S$ and one such vertex is
distinguished as a root. Decorated vertices which are not a root are
called the ends. Furthermore such a tree has a given embedding into
the plane.
\end{defn}

The embedding into the plane defines a canonical ordering on the
edges which is obtained by going counter clockwise starting with the
root edge. The root defines a direction on each edge which is the
direction away from the root.

\begin{defn}{\bf (Tree algebra $\T_S$).} The free graded commutative algebra generated by the rooted
$S$-decorated planar trees, where each tree has weight equal to its
number of edges, is denoted by ${\mathcal T}_S$. There is a
bi-grading where a rooted tree in ${\mathcal T}_S^n(r)$ has $r$ ends
and $2r-n$ edges. (Notice that the differential grading is different
with the number of edges but is congruent to it modulo 2, so for
sign purposes both give the same value).
\end{defn}
\begin{defn}{\bf (Contraction $T/e$).}
 For a rooted $S$ decorated
planar tree $T$, with a given edge $e$, the contracted tree $T/e\in
{\mathcal T}_S$ is defined as follows:
\begin{enumerate}
\item If $e$ is an internal edge then $T/e$ is obtained by
contracting the edge $e$, keeping the root and decoration and
embedding in the plane as before.
\item If $e$ is an external edge (i.e. an edge with a decorated
vertex) then $T/e$ is obtained in several stages. First remove $e$
and its two vertices. Denote the connected components by $T_1,\dots,
T_k$ (in the order dictated by the embedding), add a vertex to the
open edges of these trees and decorate it by the decoration of the
vertex of $e$. For the trees that do not have a root, make this
newly added vertex their root. $T/e$ is the product of rooted $S$
decorated trees $T_1\cdots T_k$.
\end{enumerate}
\end{defn}
 Two examples are given in the Figure 1.

\begin{figure}
\begin{center}
\epsfig{file=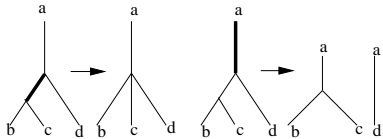, height=50pt}\caption{Two examples of
contraction of edges}\label{fig:contraction}
\end{center}
\end{figure}

\begin{defn}{\bf (Differential).}
For a rooted $S$ decorated planar tree $T$ in ${\mathcal T}_S$ the
differential is given by
\begin{equation*}
dT=\sum_{i=1}^{e(T)}(-1)^{i-1}(T/e_i).
\end{equation*}
Here the ordering is given by the embedding by going counter
clockwise starting from the root and $e(T)$ is the number of edges
of $T$. This will be extended by Leibniz rule to all of ${\mathcal
T}_S$.
\end{defn}
It is easy to check that $d^2=0$. Therefore ${\mathcal T}_S$ becomes
a commutative DGA.
\\
\begin{defn}
Let $S$ and $S'$ be two sets. The DGA of double decorated rooted
trees ${\mathcal T}_{S,S'}$ is the free graded algebra generated
by rooted $S$-decorated planar trees together with a decoration of
edges with values in $S'$. We require two conditions:
\begin{enumerate}
\item The two decoration of the leaves are distinct (we do not
make this assumption for the root).

\item For any edge decorated with $s'$ and final vertex $v$,
either $v$ is a decorated vertex or there is an edge starting from
$v$ with the same edge decoration.

\item Any path between two labeled vertex with the same label $s$,
passes through the start vertex of an edge decorated with $s$.
\end{enumerate}

Its easy to check that the previous differential defines a
differential structure on this algebra if we make the following
modification: we do not collapse leaves with vertex decoration the
same as the edge decoration of their mother edge. Notice that the
collapse of such an edge will produce a tree that does not satisfy
the condition (1) above, nevertheless it does not produce any
problem since this produces another tree with a root vertex
decorated by same label as the root edge and the cycle associated
to such a tree will be zero in the following definition.

\end{defn}

 Bearing in mind the integral $\int_a^c\frac{dt}{t-b}$,
 for an oriented edge $e$ with origin vertex labeled by $a$ and the
final vertex labeled by $b$ (they can be elements of $F$ or
variables) and the edge labeled by $c$ define the function

\begin{equation*}
 f(e)=\frac{a-b}{c-b}.
\end{equation*}
\\
 We will define the following map which generalizes a construction
of Gangl, Goncharov and Levin in \cite{GGL}:
\begin{defn}{\bf (Forrest cycling map).}
Decorate the internal vertices of elements of $\T_{S,S'}$ by
independent variables. Give each edge the direction that points
away from the root. The forrest cycling map to an element $T$ of
$\T_S^n(r)$ associates the following cycle of codimension $r$
inside $\Box_F^{2n-r}$:
$$\rho(T)=\mbox{Alt}_{2n-r}(f(e_1),\dots,f(e_{2r-n})).$$
Here $e_1,\dots,e_{2r-n}$ are the edges of $T$ with the ordering
induced from its embedding.
\end{defn}

Let us explain this in more details. Let $I$ be the set of
internal vertices. Then $(f(e_1),\dots,f(e_{2r-n}))$ is a rational
function from $$\phi:({\Bbb P}^1)^I\dashrightarrow ({\Bbb
P}^1)^{2r-n}$$. $\rho(T)$ is defined by the alternation of the
following cycle:
$$\phi_*\left(({\Bbb P^1})^I\right)\cap \Box_F^{2r-n}.$$
This is a cycle of dimension equal to the number of internal
vertices and hence the codimension equal to the number of end
vertices.

\begin{prop}
The map $\rho:\T_{S,S'}\To {\N}$ is a well-defined morphism of
DGA's.
\end{prop}
{\bf Proof.} The compatibility with product:
$\rho(T_1T_2)=\rho(T_1)\rho(T_2)$ is clear since $T_1$ and $T_2$
have different variables for their internal vertices. Therefore to
show that $\rho(T)$ is admissible and $d\rho(T)=\rho(dT)$, we can
assume that $T$ is connected. To show that $\rho(T)$ is
admissible, first notice that this cycle does not intersect with
the faces obtained by letting some coordinate equal to $\infty$.
If the function associated to the edge $e$ from $a$ to $b$ with
label $c$ is equal to $\infty$ then we are in one of the following
two cases. Either $b$ is a variable and $b=c$, in which case the
edge $e'$ starting from $b$ with label $c$, gives the value $1$
which is not in $\Box_F$ (such an edge exist because of condition
(2) above), or $a$ is a variable and $a=\infty$ , in which case
the edge $e'$ ending at $a$
 gives the value $1$. Now if we restrict the cycle $\rho(T)$ to the face with
 $i^{\mbox{th}}$-coordinate (corresponding to the edge $e_i$ from $a$ to $b$)
 equals
 zero, this means that we have to let $a=b$. Therefore we get $\rho(T/e_i)$.  By induction on
 the number of edges we see that $\rho(T)$ is admissible and as a
 side we have proved that $d\rho(T)=\rho(dT)$. Notice that for the induction to work we need to assume validity of the condition 3, this is because
 a single edge with labels $a$ and $a$ gives the cycle $\{0\}-\{\infty\}$ which is not
admissible. \qed

\begin{rem}\label{xmass}
The map $\rho$ factors through a quotient of $\T_{S,S'}$ which we
denote by $\tilde{\T}_{S,S'}$. Two decorated trees $T_1$ and $T_2$
in $\T_{S,S'}$ which can be transformed to each other by an
automorphism of trees respecting the labels have the relation
$T_1=\epsilon T_2$ where $\epsilon=1$ if the automorphism is an
even permutation of the edges and $\epsilon=-1$ otherwise.
\end{rem}

\section{Iterated Integrals and Trees.}

Let $S$ be a non-empty subset of ${\mathbb A}^1$. For $a_i\in S$,
for $i=0,\dots,n+1$ define $\tilde{t}(a_0;a_1,\dots,a_n)\in
\T_S^1(n)$ as the sum of all rooted planar 3-valent trees with $n$
leaves decorated by $a_1,\dots,a_n$ (in this order) and its root
decorated by $a_0$. We define
$$t(a_0;a_1,\dots,a_n;a_{n+1}):=\tilde{t}(a_0;a_1,\dots,a_n)-\tilde{t}(a_{n+1};a_1,\dots,a_n).$$
In this definition $t(a;b)=0$ by our convention.
\begin{prop}\label{tree}
The differential of $t(a_0;a_1,\dots,a_n;a_{n+1})$ is
given by
$$-\sum_{0\le i<j\le
n}t(a_0;a_1,\dots,a_i,a_{j+1},\dots;a_{n+1})t(a_i;a_{i+1},\dots,a_j;a_{j+1})$$
where the first term is obtained by removing $a_{i+1}$ up to $a_j$.
\end{prop}

{\bf Proof.} The contraction of internal edges will cancel each
other. Refer to the figure 2.

\begin{figure}[h]
\begin{center}
\epsfig{file=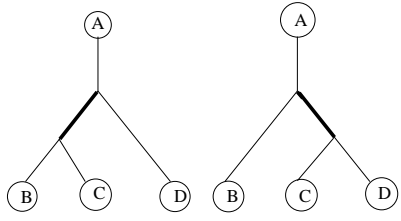, height=80pt}
\caption{The contraction of the thick edge will give the same graph but with opposite signs}\label{fig:cancel}
\end{center}
\end{figure}

In the Figure 3, let $B$ denote the decoration $a_{i+1},\dots,a_j$ for $0\le i<j\le n$ and $x=a_i$ and $y=a_{j+1}$,
 the contraction of the thicken edges with the appropriate sign will give the following term:

$$-t(a_0;\dots,a_i,a_{j+1},\dots;a_{n+1})t(a_i;a_{i+1},\dots,a_{j+1})$$

\begin{figure}[h]
\begin{center}
\epsfig{file=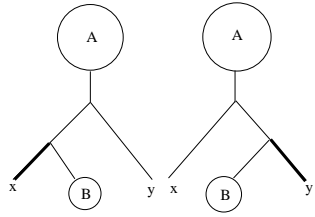, height=80pt}\caption{The contraction of the
thick leaf edges will give the desired term in the
differential}\label{fig:leaf}
\end{center}
\end{figure}
A similar argument for contraction of the root edge will prove the
theorem.\qed

\begin{defn}{\bf (Admissible decomposition).}
An admissible decomposition of $(a_0;a_1,\dots,a_n;a_{n+1})$ is an
ordered decomposition $D=P_1\cup \dots\cup P_k$ of the regular
polygon with vertices (with clockwise order) $a_0,\dots,a_{n+1}$
into sub polygons $P_i$ by diagonals. The ordering should satisfy
the following admissibility condition. If for $i<j$, $P_i$ and $P_j$
have a common edge and their union has vertices $a_{i_1},\dots,
a_{i_m}$ in clockwise ordering starting with the vertex with the
smallest index, then the edge $a_{i_1}a_{i_m}$ should belong to
$P_i$.
\\
For a sub polygon $P$ with vertices $a_{i_1},\dots, a_{i_m}$ ordered
as above we let
$$t(P):=t(a_{i_1};\dots;a_{i_m}),$$
and for a decomposition $D=P_1\cup\dots \cup P_k$ we let:
$$t(D):=[t(P_1)|\dots|t(P_k)]\in B(\T)^0.$$
\end{defn}
\begin{cor}\label{bump}
The element $T(a_0;\dots;a_{n+1}):=\sum_D t(D)\in B(\T)^0$ (where
the sum is over all the admissible decomposition of ${\bf a}$ as
above) has differential zero so defines an element in $\h^0B(\T)$
denoted again by $T({\bf a})$.
\end{cor}
{\bf Proof.} The differential of $t(D)$ for an admissible
decomposition $D=P_1\cup\dots\cup P_k$ has two part, the external
differential:
$$\sum_{i=1}^k [t(P_1)|\dots|d t(P_i)|t(P_{i+1})|\dots|t(P_k)]$$
and the internal differential
$$\sum_{i=1}^{k-1}
[t(P_1)|\dots|t(P_i)t(P_{i+1})|\dots|t(P_k)]$$ but according to
theorem 4.1:
$$dt(P_i)=-\sum_j t(P_{ij})(P_{ij}')$$
where the sum is taken over all possible division of the polygon
$P_i$ into two sub-polygon by a diagonal. This shows that the
internal differential for admissible decompositions of order $k$
will be canceled by the external differential of decompositions of
order $k-1$. Therefore the differential of the sum over all
possible decompositions vanishes.\qed

Let $\T_S'$ be the sub-DGA of $\T_S$ generated by admissible
trees, i.e. those trees that the non-zero labels do not repeat.
There is a map of DGA's :
$$dec: \T_S'\To \T_{S,\{0,1\}}$$
given by decorating the leaves with zero label by 1 and the rest
of edges by 0.

\begin{defn} For generic sequence $(a_0;a_1,\dots,a_n;a_{n+1})$, i.e. a sequence with non repeating non-zero terms,
define the motivic analogue of the iterated integral
$$(-1)^n\int_{a_0}^{a_{n+1}}\frac{dt}{t-a_1}\circ
\dots\circ\frac{dt}{t-a_n}$$ by
$$\rho(dec(T(a_0;a_1,\dots,a_n;a_{n+1}))).$$
This is an element of $\chi_F(n)=\h^0B(\N)(n)$ and will be denoted
by ${\mathbb I}(a_0;\dots;a_{n+1})$. Here
$\rho:\widetilde{\T}_{S,\{0,1\}}\To\N$ and it induces a morphism
$\rho:\h^0B(\widetilde{\T}_{S,\{0,1\}})\To\h^0B(\N)$.
\end{defn}

\section{Proof of the Theorems \ref{properties} and \ref{coproduct}.}

In this section we will prove theorems \ref{properties} and
\ref{coproduct} from the introduction. It shows that
$\I(a_0;\dots;a_{n+1})$ deserves to be called motivic iterated
integral.
\\
In fact we have constructed $\I(a_0;\dots;a_{n+1})$ as the image
(under the morphism $\rho\circ dec$ of an element
$T(a_0;\dots;a_{n+1})$ in $B(\tilde{\T})^0$ with zero
differential. We will prove the identities of theorem
\ref{properties} as identities in $B(\tilde{\T})^0$ with $\I$
replaced by $T$. This obviously implies the corresponding
identities in $\h^0B(\N)$.
\\
 We take $\I(a;b)=1$ as definition. Since $t(a;\dots;b)=0$ for $a=b$ it
 follows from the construction that $T(a;\dots;b)=0$ for $a=b$. This
 proves the first part of the theorem.
 \\
We now prove the shuffle relation. Let $D=P_1\cup\dots\cup P_k$ be
an admissible decomposition of $(a;s_1,\dots,s_{n+m};b)$. For a
permutation $\sigma$ of $s_1,\dots,s_{n+m}$ we denote by $\sigma(D)$
the decomposition of $$(a;s_{\sigma(1)},\dots,s_{\sigma(n+m)};b)$$
obtained by applying $\sigma$ to the vertices of the polygons, with
the ordering remained as before. Suppose that one of the polygons
$P_i$ of the decomposition $D$ has vertices (in order)
$s_{i_0},\dots,s_{i_{l+1}}$ such that the sequence $i_1,\dots,i_l$
is mixed i.e. it has both numbers less than or equal to $n$ and
numbers bigger than $n$. we show that
$$\sum [t(P_1)|\cdots|t(\sigma(P_i))|\cdots|t(P_k)]=0$$
where the sum is taken over all shuffles of $(i_1,\dots,i_l)$ for
the indices less than or equal to $n$ and those bigger than $n$. To
prove this for any 3-valent tree with decoration
$s_{i_1},\dots,s_{i_l}$ (with this order or an order obtained by a
shuffle of them) we define a dual tree that cancel it. The
separating edge of a tree decorated by $s_{i_1},\dots,s_{i_l}$ is
the first edge (using the ordering of the edges) that has the
property one of the  sub-trees that grows out of this edge has
decoration by indices less than or equal to $n$ and the other one
has decoration by indices bigger than $n$. For a tree $T$ with the
separating edge $e$ we define its dual the tree obtained from this
tree by switching the two sub-trees of $e$, this is obtained by a
shuffle on the decorating and clearly cancels $T$, in the modified
algebra $\tilde{T}$. Therefore we have to only consider the
decompositions of $(a;s_{\sigma(1)},\dots,s_{\sigma(n+m)};b)$ for
$(n,m)$-shuffles $\sigma$, that are clean, i.e. each sub division
polygon with vertices $s_{i_0},\dots,s_{i_{l+1}}$ has the property
that all the indices $i_1,\dots,i_l$ are either bigger than $n$ or
less than or equal to $n$. Its easy to see that the sum of $T(D)$
over this decompositions is equal to
$\T(a;s_1,\dots,s_n;b)T(a;s_{n+1},\dots,s_{n+m};b)$. This finishes
the proof of the shuffle relation.
\\
We now prove the path composition formula. For an admissible
decomposition $D=P_1\cup \cdots\cup P_m$ of $(a;s_1,\dots,s_n;b)$,
let $P_1=P_{i_0},P_{i_1},\dots,P_{i_l}$ be the sub-polygons with $b$
as a vertex. Using the trivial identity:
$$t(s_{j_1};\dots,s_{j_p};b)=t(s_{j_1};\dots,s_{j_p};c)+t(c;\dots,s_{j_p};b)$$
we can replace $t(P_{i_p})$ as the sum of two terms
$t(P_{i_p}')+t(P_{i_p}'')$ where $P_{i_p}'$ is obtained by replacing
the last vertex $b$ by $c$ in $P_{i_p}$ and similarly $P_{i_p}''$ is
obtained by replacing the first vertex (the one with the smallest
index) by $c$. Therefore we can write $t(D)$ as:
$$\sum_{k=0}^l
[t(P_1')|\dots|t(P_{i_{k-1}}')|\dots|t(P_{i_k}'')|t(P_{i_k+1})|\dots|t(P_m)].$$
The $k^{th}$ term of this sum when $P_{i_k}=(c;s_i,\dots;b)$ belongs
to the product: $$\I(a;s_1,\dots,s_{i-1};c)\I(c;s_i,\dots;b).$$ It
is clear that this way we get an identification between the terms in
$\I(a;s_1,\dots,s_n;b)$ and
$$\sum_{i=1}^{n+1}\I(a;s_1,\dots,s_{i-1};c)\I(c;s_i,\dots,s_n;b).$$ This finishes
the proof of path composition.
\\
The proof of coproduct formula goes as follows. Recall the
coproduct:
$$\Delta [a_1|\dots|a_n]=\sum_{r=0}^n
[a_1|\dots|a_r]\otimes[a_{r+1}|\dots|a_n].$$  Thus for a
decomposition $D=P_1\cup\cdots\cup P_m$ we have:

$$\Delta t(D)=\sum_{r=0}^m [t(P_1)|\dots|t(P_r)]\otimes
[t(P_{r+1})|\dots|t(P_m)].$$ Taking the sum over all decompositions
such that $P_1\cup\cdots\cup P_r$ is the polygon
$(a;s_{i_1},\dots,s_{i_k},b)$ for a fixed sequence
$0=i_0<i_1<\cdots<i_k<i_{k+1}=n+1$ and varying $r$ we get the term

$$T(a;s_{i_1},s_{i_2},\dots,s_{i_k};b)\otimes \prod
_{j=0}^k T(s_{i_j};s_{i_j+1},\dots,s_{i_{j+1}-1};s_{i_{j+1}})
$$
in the coproduct. This finishes the proof of coproduct formula. To
prove the antipode relation note that the mirror dual of a tree with
$n$ ends is equal to that tree time $(-1)^{n-1}$ in the modified
algebra $\tilde{T}$, so if we switch $a$ and $b$ as well a sign of
$(-1)^n$ will appear. This finishes proof of the theorem
\ref{properties} and \ref{coproduct}.\qed

\section{Review of the Hodge Realization.}
 In this section we review the construction of Bloch-K\v ri\v z for
associating a framed MHTS to an element of $\h^0B(\N)$. In fact for
our application it is only enough to recall \S 8 of \cite{BK} where
they give the Hodge realization for $\h^0B(\N')$ for a particular
sub-DGA $\N'$ of $\N$.
\\
Let $\omega(n,r)$ be a collection of real oriented subvarieties of
$({\Bbb P}^1)^n$ of codimension $2r$. Assume that for
$S_1\in\omega(n,r)$ and $S_2\in \omega(m,s)$, $S_1\times
S_2\in\omega(n+m,r+s)$, and for $\sigma\in G_n=(\Z/2)^n\rtimes
\Sigma_n$, $\sigma(S_1)\in\omega(n,r)$. Moreover assume that the
intersection of an element of $\omega(n,r)$ with a hyperplane
$t_i=0$ or $\infty$ belongs to $\omega(n-1,r)$.
\\
Define a DGA with an Adams grading by:

$${\mathcal D}^n(r):=\mbox{Alt}\varinjlim\h_{2r-2n}(S\cup{\mathcal
J}^{2r-n},{\mathcal J}^{2r-n})$$

where ${\mathcal J}^n$ is the union of all the codimension 1 hyper
planes  of $({\Bbb P}^1)^n$ obtained by letting one coordinate equal
to 1. The limit is taken over all $S\in\omega(2r-n,r)$. This has a
natural structure of DGA (refer to \cite{BK}, \S 8).
\\

\begin{defn}\label{admissible} The pair $(\N',{\mathcal D})$ as above, is called admissible if
the following conditions hold:
\begin{enumerate}
\item The Adams graded pieces $\h^0B(\N')(r)$ are finite dimensional
$\Q$-vector spaces.
\item All the support subvarieties of the elements of $\N'$ belong
to $\omega(*,*)$. Hence there is a natural morphism of DGA's:
$$\sigma:\N'\To {\mathcal D}$$
given by the fundamental class.
\item The map $\lambda_0: {\mathcal D}\To \C[x]$ given by sending an
element $c\in {\mathcal D}^0(n)$ to
$$(\int_c \omega_{2n})\cdot x^n$$ and ${\mathcal D}^i$ for $i\ne 0$ to zero,
is a well-defined morphism of DGA's. Here $\omega_n=(2\pi
i)^{-n}\frac{dz_1}{z_1}\wedge \dots\wedge \frac{dz_{n}}{z_{n}}$ and
$\C[x]$ is a DGA with trivial differential and concentrated in
degree zero with an Adams grading given by powers of $x$.
\item If $\sum_i [a_1^i|\dots|a_{r_i}^i]\in B(\N')^0$ has zero differential then the element $$\sum_i \sigma(a_1^i)[a_2^i|\dots|a_{r_i}^i]\in
B({\mathcal D},\N')^1$$ is in $d(B({\mathcal D},\N')^0)$.
\item The map $\alpha$ defined in \cite{BK} (8.11)(There is a missprint: it is written $\h^0B({\mathcal D}/\N'(r))$, instead of $\h^0({\mathcal D}/\N'(r))$):
$$\alpha:\h^0({\mathcal D}/\N'(r))\To \h_{2r}(({\Bbb
P}^1-\{0,\infty\})^{2r}\cup N(r)\cup{\mathcal J}^{2r},{\mathcal
J}^{2r})/\mbox{Im}\gamma$$ preserves integration of $\omega_{2r}$.
\\
Here $N(r)$ is the union of codimension $r$ algebraic subvarieties
of $\Box_{\C}^{2r}$ that intersect all the faces properly and
$\gamma$ is the map:
$$\gamma: \oplus_H\h_{2r}({\mathcal J}^{2r},H)\To \h_{2r}(({\Bbb
P}^1-\{0,\infty\})^{2r}\cup N(r)\cup{\mathcal J}^{2r},{\mathcal
J}^{2r})$$ induced by inclusion, and $H$ runs over hyper planes with
one coordinate equal to one.
\end{enumerate}
\end{defn}
Let us quickly recall some basic notions from mixed Hodge-Tate
structures.
\begin{defn} A MHTS $H$ is a finite dimensional $2\Z$-graded $\Q$-vector
space $H_{dR}=\oplus_n H_{2n}$ together with a $\Q$-subspace $H_B$
of $H_{dR}\otimes \C$ such that for all $m$:
$$\mbox{Im}\left(H_B\cap (\bigoplus_{n\le m} H_{2n}\otimes \C)\stackrel{Proj}\To
H_{2m}\otimes \C\right)= H_{2m}.$$ A morphism $f$ is a graded
morphism $f_{dR}$ from $H_{dR}$ to $H_{dR}'$ such that
$f_{dR}\otimes 1:H_{dR}\otimes \C\To H_{dR}'\otimes \C$ sends $H_B$
to $H_B'$, this induced map is denoted by $f_B$.
\end{defn}
The notion of framed MHTS was introduced by Beilinson , Goncharov,
Schechtman and Varchenko in \cite{BGSV}. This will give a concrete
way of thinking about the Hopf algebra of the coordinate ring of the
tannakian Galois group. \\
\begin{defn}
 A framing on a MHTS $H$ is the choice of a frame vector $v\in H_{2r}$ and a co-frame vector $\hat{v}\in \mbox{Hom}(H_{2s},\Q)$.
 Two framed MHTS $(H,v ,\hat{v})$ and $(H',v',\hat{v}')$ such that for an integer $n$ there is a morphism from $H(n)\To H'$ that
 respects the frames, are called equivalent. This relation generates an equivalence relation. The equivalence classes of framed MHTS
 is denoted by $\chi_{MHTS}$. This is graded where $[(H,v,\hat{v}]$ where $v\in H_{2r}$ and $\hat{v}\in \mbox{Hom}(H_{2s},\Q)$ has degree $r-s$.
 It is easily shown that the negative weights vanish. The following definitions of product and coproduct make $\chi_{MHTS}$ into a Hopf algebra:
$$[(H,v,\hat{v})]\cdot [(H',v',\hat{v}')]=[(H\otimes H',v\otimes v',\hat{v}\otimes \hat{v}')]$$
$$\Delta=\bigoplus_m \Delta_m :\;\;\;  \Delta_m[(H,v,\hat{v})]=\sum_i [(H,v,\hat{b}_i)]\otimes [(H,b_i,\hat{v}')]$$
where $b_i$ is a basis for $H_{2m}$ and $\hat{b}_i$ is the dual
basis for $\mbox{Hom}(H_{2m},\Q)$. It can be shown that the category
of finite dimensional graded comodules over $\chi_{MHTS}$ is
equivalent to the category of MHTS's.
\end{defn}
Under the conditions explained above on ${\mathcal D}$  Bloch-K\v
ri\v z prove that:
\begin{thm}\label{Hodge}
(\cite{BK}) The graded $\Q$-vector space $\oplus_{0\le r\le n}
\h^0B(\N')(r)$ where $\h^0B(\N')(r)$ has degree $2r$, together with
the subspace defined by the image of the map $\lambda$ defined as
composition
$$\h^0B({\mathcal
D},\N')\xrightarrow{B(\lambda_0,id)}
\h^0B(\C[x],\N')\stackrel{x=1}\To \h^0B(\N')_{\C}$$ is a MHTS. We
denote it by $H(\N',n)$.
\end{thm}
\begin{defn}\label{real}
The realization map $\mbox{Real}'_{MHTS}: \h^0B(\N')(n)\to
\chi_{MHTS}(n)$ is given by sending $a$ to the framed MHTS
$H:=H(\N',n)$ framed by $$a\in \h^0B(\N')(n)=H_{2n}$$ and the
augmentation isomorphism $$\h^0B(\N')(0)\rightarrow \Q$$ which is an
element of $\mbox{Hom}(H_{0},\Q)$.
\end{defn}
It is proven in \cite{BK} that:
\begin{thm}
The realization map $\mbox{Real}'_{MHTS}$ is morphism of graded Hopf
algebras. It is equal to the composition of $\h^0B(\N')\To
\h^0B(\N)$ with the general Hodge realization
$\mbox{Real}_{MHTS}:\h^0B(\N)\To \chi_{MHTS}$ constructed in \S 7 of
\cite{BK}.
\end{thm}

\section{Hodge Realization of the Motivic Iterated Integrals.}

In this section we will calculate the Hodge realization of
${\mathbb I}({\bf a})\in \h^0B(\N')$, where ${\bf
a}=(a_0;a_1,\dots,a_n;a_{n+1})$ is a {\it{generic}} sequence in
$S$. We will do this parallel to the special case
$${\mathbb I}_n(a):={\mathbb I}(0;\underbrace{1,0,\dots,0}_n;a)$$
 for
$a\in F$ which is done in \cite{BK}. Here $\N'=\N_{\bf a}$ is the
sub-DGA of $\N$ generated by the cycles $\rho\circ dec\circ t({\bf
a}')\in \N^1(n)$ for all sub-sequences ${\bf a}'$
of ${\mathbf a}$. \\
Notice that by definition $\N'$ is a minimal DGA, i.e. it is
connected and $d(\N')\subseteq {\N'}^+\cdot {\N'}^+$.
\begin{lem}\label{generate}
The graded Hopf algebra $\h^0B(\N')$ is generated as a DGA by the
elements ${\mathbb I}({\bf a}')$ for all sub-sequences of ${\bf a}$.
\end{lem}
{\bf Proof.} It follows from Theorem 6.3 in \cite{KM}.\qed
\\
 Once and for all choose a path $\gamma$ with
interior in ${\mathbb A}^1(\C)-S$ from the tangential base point
$a_0$ to the tangential base point $a_{n+1}$. For $1\le i\le n$ we
define intermediate cycles $\eta_i({\bf a})$ inside
$\Box_{\C}^{2n-i}$ with (real) dimension $2n-i$. It is defined as
follows. Consider the sum of all rooted planar forests with $i$
connected components and $n$ leaves, decorated by $a_1,\dots,a_n$.
Its roots are decorated by variables
$\gamma(s_1),\dots,\gamma(s_i)$ for $0\le s_1\le\dots\le s_i\le
1$. Apply the morphism $\rho\circ dec$ to this which will give an
oriented (real) cycle (with boundary) of dimension $2n-i$ inside
$\Box_{\C}^{2n-i}$. Let $\Gamma$ be a small disk around zero in
$\Box_{\C}$ with its canonical orientation. We define:
$$\tau_i({\bf a})=(\delta\eta_i({\bf a}))\cdot \Gamma\cdot (\delta\Gamma)^{i-1}+(-1)^i\eta_i({\bf a})\cdot(\delta\Gamma)^i$$
where $\cdot$ denotes the usual alternating product and $\delta$
denotes the topological boundary defined for a cycle
$f(s_1,\dots,s_n)$ by:

\begin{align*}
(\delta f)(s_1,\dots
s_{n-1})=-f(0,s_1,s_2,\dots,s_{n-1})+f(s_1,s_1,s_2,\dots,s_{n-1})-\cdots
\\
 \cdots+(-1)^n
f(s_1,\dots,s_{n-1},s_{n-1})+(-1)^{n+1}f(s_1,\dots,s_{n-1},1).
\end{align*}
Finally denote
$$\xi_{\gamma}({\bf a})=\sum_{k=1}^n \tau_k({\bf a}).$$
Let $\omega(*,*)$ be the subvarieties which support the cycles
$\xi({\bf a}')\rho({\bf a}'')$ and $\rho({\bf a}')$ for all sub
sequences of ${\bf a}'$ of ${\bf a}$, here ${\bf a}''$ is the
complement of ${\bf a}'$ inside $\bf a$. Denote by ${\mathcal
D}={\mathcal D}_{\bf a}$ the corresponding DGA obtained from
$\omega(*,*)$. Hence $\xi_{\gamma}({\bf a)})\in {\mathcal D}^0(n)$.
\\
\begin{rem}
>From now on we will use $\rho$ instead of $\rho\circ dec\circ t$
for the ease in typing.
\end{rem}
 The crucial technical result about the cycle
$\xi_{\gamma}({\bf a})$ is
\begin{prop}\label{diff}

The differential of $\xi_{\gamma}({\bf a})$ is given by the
formula
\begin{equation*}
\rho({\bf a})-\sum_{0\le i<j\le
n}\xi_{\gamma}(a_0;\dots,a_i,a_{j+1},\dots;a_{n+1})\rho(a_i;\dots;a_{j+1}).
\end{equation*}
\end{prop}
As in corollary \ref{bump} for an admissible decomposition
$D=P_1\cup \dots\cup P_k$ let
$\xi_{\gamma}(D):=\xi_{\gamma}(P_1)[\rho(P_2)|\dots|\rho(P_k)]\in
B({\mathcal D},\N')$.
\begin{cor}\label{Z} Define $Z_{\gamma}({\bf a}):=1\cdot [\I({\bf a})]+\sum_D \xi_{\gamma}(D)$ where
the sum is over all admissible decompositions of ${\bf a}$. Then
$d(Z_{\gamma}({\bf a}))=0$. Therefore $Z_{\gamma}({\bf a})$ defines
an element of $\h^0B({\mathcal D},\N')$.
\end{cor}
{\bf Proof.} Note that:
$$d(1\cdot[\I({\bf a})])=-\sum_{D}
\rho(P_1)[\rho(P_2)|\dots|\rho(P_k)]$$ where the sum is taken over
all admissible decompositions $D=P_1\cup\dots\cup P_k$ of $\bf a$.
On the other hand for $D$ as above:
\begin{align*}
d_{ext}\xi_{\gamma}(D)=&(d\xi_{\gamma}(P_1))[\rho(P_2)|\dots|\rho(P_k)]+\sum_{i=2}^k
\xi_{\gamma}(P_1)[\rho(P_2)|\dots|d\rho(P_i)|\dots|\rho(P_k)]
\\
=& \rho(P_1)[\rho(P_2)|\dots|\rho(P_k)]-\sum_j
\xi_{\gamma}(P_{1j})\rho(P_{1j}')[\rho(P_2)|\dots|\rho(P_k)]
\\
&+\sum_{i=2}^k
\xi_{\gamma}(P_1)[\rho(P_2)|\dots|d\rho(P_i)|\dots|\rho(P_k)]
\end{align*}
Also
\begin{align*}
d_{int}\xi_{\gamma}(D)=&\xi_{\gamma}(P_1)\rho(P_2)[\rho(P_3)|\dots|\rho(P_k)]
\\
&+\sum_{i=2}^{k-1}
\xi_{\gamma}(P_1)[\rho(P_2)|\dots|\rho(P_i)\rho(P_{i+1})|\dots|\rho(P_k)].
\end{align*}

Since $d(\rho(P_i))=-\sum_j\rho(P_{ij})\rho(P_{ij}')$ for all
divisions of $P_i$ into two admissible polygons $P_{ij}$ and
$P_{ij}'$. It will be seen that the third term in the external
differential for $k$ will be cancelled by the second term of the
internal differential for $k+1$. The second term of the external
differential for $k$ will be cancelled by the first term of the
internal differential for $k+1$. Finally the first term will be
cancelled by the differential of $1\cdot [\I(\bf a)]$.\qed
\\
\\
{\bf Proof of Proposition \ref{diff}}. For the proof we need the
following two lemmas:
\begin{lem}\label{lem1}
The differential of $\eta_k({\bf a})$ is given by:
$$\delta \eta_{k+1}({\bf a})-(-1)^{k}\sum_{0\le i< j\le
n}\eta_k(a_0;a_1,\dots,a_i,a_{j+1},\dots;a_{n+1})\rho(a_i;\dots;a_{j+1}).$$
\end{lem}
\begin{lem}\label{lem2}
The differential of $\delta\eta_k({\bf a})$ is given by:
$$(-1)^{k-1}\sum_{0\le i<j\le
n}\delta\eta_k(a_0;\dots,a_i,a_{j+1},\dots;a_{n+1})\rho(a_i;\dots;a_{j+1}).$$
\end{lem}
Assuming these lemmas, it follows:
\begin{align*}
d\tau_k({\bf a})=& d\delta \eta_k({\bf a})\cdot \Gamma\cdot
\delta\Gamma^{k-1}+(-1)^k\delta\eta_k({\bf a})\cdot
\delta\Gamma^{k-1}+(-1)^k d\eta_k({\bf a})\cdot \delta\Gamma^k
\\
=& -\sum_{0\le i<j\le
n}\Bigl\{\delta\eta_k(a_0;\dots,a_i,a_{j+1},\dots;a_{n+1})\cdot
\Gamma\cdot\delta\Gamma^{k-1}
\\
&+(-1)^k\eta_k(a_0;\dots,a_i,a_{j+1},\dots;a_{n+1})\cdot\delta\Gamma^k\Bigr\}\cdot\rho(a_i;\dots;a_{j+1})
\\
&+(-1)^k\delta\eta_k({\bf a})\cdot
\delta\Gamma^{k-1}+(-1)^k\delta\eta_{k+1}({\bf a})\cdot
\delta\Gamma^{k}
\\
=&-\sum_{0\le i<j\le n}\tau_k
(a_0;\dots,a_i,a_{j+1},\dots;a_{n+1})\cdot\rho(a_i;\dots;a_{j+1})
\\
&+(-1)^k\delta\eta_k({\bf a})\cdot
\delta\Gamma^{k-1}+(-1)^k\delta\eta_{k+1}({\bf a})\cdot
\delta\Gamma^{k}.
\end{align*}
Adding for $k=1,\dots,n$ we will get:
$$d\xi_{\gamma}({\bf a})=-\delta\eta_1({\bf a})-\sum_{0\le i<j\le
n}\xi_{\gamma}
(a_0;\dots,a_i,a_{j+1},\dots;a_{n+1})\cdot\rho(a_i;\dots;a_{j+1})$$
but since $\delta\eta_1({\bf a})=-\rho({\bf a})$ this finishes the
proof of theorem \ref{diff} modulo the two lemmas.\qed
\\
{\bf Proof of lemma \ref{lem1}.} To calculate differential of
$\eta_k({\bf a})$ we have to calculate the differential of the sum
over all 3-valent forests with $k$ components decorated by
$a_1,\dots,a_n$ for the ends and $\gamma(s_1),\dots,\gamma(s_k)$
for the roots, where $s_i$'s are variables so that $0\le
s_1\le\dots\le s_k\le 1$. To calculate this differential we have
to consider the contraction of the edges. As in the proof of
proposition \ref{tree} the only edges that contribute are external
edges. The contraction of the roots after applying the morphism
$\rho:\T\To\N$ will give all the terms in $ \delta\eta_{k+1}({\bf
a})$, except the terms corresponding to $s_1=0$, i.e.
$\gamma(s_1)=a_0$ and $s_k=1$, i.e. $\gamma(s_k)=a_{n+1}$. The
contraction of the leaves can be grouped exactly like figure 3 in
the proof of proposition \ref{tree}, and this together with two
missing terms in $\delta\eta_k({\bf a})$ will give (after applying
the morphism $\rho:\T\To \N$):
$$(-1)^{k-1}\sum_{0\le i<j\le
n}\eta_k(a_0;\dots,a_i,a_{j+1},\dots;a_{n+1})\cdot
\rho(a_i;\dots;a_{j+1}).$$ This finishes the proof of lemma
\ref{lem1}.\qed
\\
{\bf Proof of lemma \ref{lem2}.} From lemma \ref{lem1} it follows
that
\begin{align*}
d\delta\eta_{k+1}({\bf a})=&(-1)^{k}d\left(\sum_{0\le i<j\le
n}\eta_k(a_0;\dots,a_i,a_{j+1},\dots;a_{n+1})\cdot
\rho(a_i;\dots;a_{j+1})\right)
\\
=&(-1)^{k}\sum_{0\le i<j\le
n}\Bigl\{\delta\eta_{k+1}(a_0;\dots,a_i,a_{j+1},\dots;a_{n+1})\cdot
\rho(a_i;\dots;a_{j+1})
\\
&+(-1)^{k+1}
d'\Bigl(\eta_k(a_0;\dots,a_i,a_{j+1},\dots;a_{n+1})\cdot
\rho(a_i;\dots;a_{j+1})\Bigr)\Bigr\}.
\end{align*}
Here $d'$ is defined as the differential on $\rho$ but its defined
by:
$$d'\eta_k(a_0;\dots;a_{m+1}):=\sum_{0\le i<j\le
m}\eta_k(a_0;\dots,a_i,a_{j+1},\dots;a_{m+1})\cdot
\rho(a_i;\dots;a_{j+1}).$$ Now it follows from anti-commutativity
that:
$$\sum_{0\le i<j\le n}d'\Bigl(\eta_k(a_0;\dots,a_i,a_{j+1},\dots;a_{n+1})\cdot
\rho(a_i;\dots;a_{j+1})\Bigr)=0,$$ so lemma \ref{lem2} is
proved.\qed
\\
\begin{lem}\label{lambda}
Assume $a_0\ne a_1$ and $a_n\ne a_{n+1}$. The image of
$\xi_{\gamma}({\bf a})\in {\mathcal D}^0(n)$ under the map
$\lambda_0$ defined in definition \ref{admissible} is given by the
iterated integral:
$$\left(-\frac{x}{2\pi i}\right)^{n} \int_{\gamma}\frac{dt}{t-a_1}\circ\dots\circ \frac{dt}{t-a_n}.$$
Furthermore  $\lambda Z_{\gamma}({\bf a})\in\h^0B(\N'_{\C})$
($\lambda$ was defined in theorem \ref{Hodge}) is:
$$\sum (-2\pi
i)^{-k}\int_{\gamma}\frac{dt}{t-a_{i_1}}\circ\dots\frac{dt}{t-a_{i_k}}\prod_{j=0}^k
\I(a_{i_j};\dots;a_{i_{j+1}}),$$ where the sum is taken over all
indices $0=i_0<i_1<\dots<i_k<i_{k+1}=n+1$ with $k=0,1,2,\dots$.

\end{lem}
{\bf Proof.} For reasons of type, the only term in
$\sum_{k=1}^n\tau_k({\bf a})$ which contributes in
$\lambda_0(\xi_{\gamma}({\bf a}))$ is $\tau_n({\bf a})$ which is
$$\tau_n({\bf a})=(\delta\eta_n({\bf a}))\cdot \Gamma\cdot (\delta\Gamma)^{n-1}+(-1)^n\eta_n({\bf
a})\cdot(\delta\Gamma)^n.$$ The first part of this sum gives zero
integral for reasons of type. Note that $\eta_n({\bf a})$ is
$$\mbox{Alt}(f(s_1,a_1),\dots,f(s_n,a_n)),$$
where
$$f(s_i,a_i)=\begin{cases}
1-\frac{\gamma(s_i)}{a_i}&\text{if $a_i\ne 0$,}
\\
\gamma(s_i)&\text{if $a_i=0$}
\end{cases}.
$$
Therefore:
$$\int_{\eta_n({\bf
a})}\frac{dz_1}{z_1}\wedge\cdots\wedge\frac{dz_n}{z_n}=\int_{\gamma}\frac{dt}{t-a_1}\circ\dots\circ\frac{dt}{t-a_n}.$$
Note that this integral is convergent due to the assumption of the
lemma. This shows that
$$x^n(2\pi i)^{-2n}\int_{\tau_n({\bf
a})}\frac{dz_1}{z_1}\wedge\cdots\wedge \frac{dz_{2n}}{z_{2n}}=
\left(-\frac{x}{2\pi i}\right)^n
\int_{\gamma}\frac{dt}{t-a_1}\circ\dots\circ\frac{dt}{t-a_n}.
$$
This proves the first part of the lemma. To prove the second part,
note that if in the set of all admissible decompositions one fixes
$P_1=(a_0,a_{i_1},\dots,a_{i_k},a_{n+1})$ then
$$\sum_{D}[\rho(P_2)|\dots|\rho(P_k)]=\prod_{j=0}^k
\I(a_{i_j};\dots;a_{i_{j+1}}).$$ This together with the first part
of the lemma, letting $x=1$ finishes the proof.\qed
\\
\begin{lem}\label{conditions}
Let ${\bf a}=(a_0;a_1,\dots;a_{n+1})$ be a generic sequence in $S$
such that $a_0\ne a_1$ and $a_n\ne a_{n+1}$. The pair $(\N_{\bf
a},{\mathcal D}_{\bf a})$ defined above is admissible in the sense
of the definition \ref{admissible}.
\end{lem}
{\bf Proof.} Since $\N'$ is finitely generated, property (1) is
clear. Property (2) follows from the definition of $\mathcal D$. To
prove property (4) we can only consider elements $\I({\bf a}')$
according to the lemma \ref{generate}, and for these elements
corollary \ref{Z} implies (4). Properties (3) and (5) follow from
the same techniques as \cite{BK}, using a Stokes formula and the
assumption on the convergence of the iterated integral, namely
$a_0\ne a_1$ and $a_n\ne a_{n+1}$.\qed

 The MHTS on the pro-unipotent fundamental torsor
$\Pi^{uni}({\Bbb A}^1-S;a_0,a_{n+1})$ is given by:
\begin{align*}
\Pi_{dR}({\Bbb A}^1-S;a_0,a_{n+1}):=&\Q\langle\langle
X_s\rangle\rangle_{s\in S}\;\;\;\mbox{weight}(X_s):=-2,
\\
\Pi_{B}({\Bbb A}-S;a_0;a_{n+1}):=&\mbox{Im}\left(\Pi^{uni}({\Bbb
A}^1-S;a_0,a_{n+1})\stackrel{\Phi}{\To}\C\langle\langle
X_s\rangle\rangle_{s\in S}\right)
\end{align*}
Here the map $\Phi$ is defined by
$$\gamma\mapsto \sum (2\pi
i)^{-n}\left(\int_{\gamma}\frac{dt}{t-s_1}\circ\cdots\circ\frac{dt}{t-s_n}\right)X_{s_1}\cdots
X_{s_n}$$ and the sum is taken over $n=0,1,\dots$ and $s_i\in S$.
\\
The Hodge analogue of the iterated integral:
$$(2\pi
i)^{-n}\int_{\gamma}\frac{dt}{t-a_1}\circ\cdots\circ\frac{dt}{t-a_n}$$
is the framed MHTS given above for $\Pi:=\Pi^{uni}({\Bbb
A}^1-S;a_0,a_{n+1})$, together with frames:
\begin{align*}
1\in &\Pi_0,
\\
(X_{a_1}\cdots X_{a_n})'\in & \mbox{Hom}(\Pi_{-2n},\Q).
\end{align*}
Here $\{(X_{s_1}\cdots X_{s_m})'\}_{s_i\in S}$  is the dual basis.
We denote this by $\I^{\mathcal H}(a_0;a_1,\dots,a_n;a_{n+1})$.
\begin{thm}\label{compare}
The Hodge realization of ${\mathbb I}({\bf a})$ is $(-1)^n
\I^{\mathcal H}({\bf a})$. Here ${\bf
a}=(a_0;a_1,\dots,a_n;a_{n+1})$.
\end{thm}
{\bf Proof.} First assume that $a_0\ne a_1$ and $a_n\ne a_{n+1}$. By
lemma \ref{conditions} $(\N',{\mathcal D})$ is admissible, hence we
can use the realization map $\mbox{Real}'_{MHTS}$ in definition
\ref{real}. Let $(J_{\bf a},\epsilon, {\mathbb I}({\bf
a})):=\mbox{Real}'_{MHTS}({\mathbb I}({\bf a}))$. We define a map
$\Pi(n) \To J_{\bf a}$ by:
$$X_{b_1}\cdots X_{b_k}\mapsto
(-1)^k\sum \prod_{j=0}^k{\mathbb
I}(a_{i_j};a_{i_j+1},\dots;a_{i_{j+1}})
$$
where the sum is taken over all indices
$0=i_0<i_1<\dots<i_k<i_{k+1}=n+1$ such that $a_{i_j}=b_j$. This
obviously extends to a graded linear map from $\Pi(n)\To J_{\bf a}$
which respects the corresponding frames. We need to show that for a
path $\gamma$ from $a_0$ to $a_{n+1}$:
$$\sum (-2\pi
i)^{-k}\int_{\gamma}\frac{dt}{t-a_{i_1}}\circ\dots\frac{dt}{t-a_{i_k}}\prod_{j=0}^k
\I(a_{i_j};\dots;a_{i_{j+1}}),$$ where the sum is taken over all
indices $0=i_0<i_1<\dots<i_k<i_{k+1}=n+1$ with $k=0,1,2,\dots$,
belongs to the Betti space of $J_a$. But this element is the image
of $Z_{\gamma}({\bf a})$ under the map $\lambda$ according to the
lemma \ref{lambda}. This completes the proof of the equivalence of
the two framed MHTS's.

We now deal with the divergent iterated integrals. Note that by our
construction $\I(0;0;1)=0$, therefore $\I(a_0;a_0;a_1)$ and
$\I(a_0;a_1;a_1)$ have correct Hodge realization corresponding to
the tangential base $\frac{\partial}{\partial t}$ for the canonical
coordinate $t$ of ${\Bbb A}^1$. Therefore using the shuffle relation
for
$$\I(a_0;\underbrace{a_0,\dots,a_0}_p;a_{n+1})\cdot
\I(a_0;a_1,\dots,a_n;a_{n+1})\cdot\I(a_0;\underbrace{a_{n+1},\dots,a_{n+1}}_q;a_{n+1})$$
where $a_0\ne a_1$ and $a_n\ne a_{n+1}$, and induction on $p+q$ it
follows that the divergent iterated integral
$$\I(a_0;\underbrace{a_0,\dots,a_0}_p,a_1,\dots,a_n,\underbrace{a_{n+1},\dots,a_{n+1}}_q;a_{n+1})$$
can be written as a linear combination of convergent iterated
integrals. Since we have a same story in the Hodge side, with exact
same formulas, it follows that the Hodge realization of the framed
MTM $\I({\bf a})$ is the framed MHTS that classically is associate
to the iterated integral:
$$(-1)^n\int_{a_0}^{a_{n+1}}\frac{dt}{t-a_1}\circ\dots\circ\frac{dt}{t-a_n}.$$
\qed

\section{Combinatorial Fundamental Group of Punctured Affine Line}

\begin{defn}

A combinatorial pro-MTM over $F$ is a graded $\Q$-vector space $V$
together with a graded co-action
$$V\To V\otimes H^0B(\tilde{\T}_F),$$
that respects the coproduct and co-unit.
\end{defn}

Let us first give the MHTS of $\Pi^{uni}({\Bbb A}^1-S;a,b)$ as a
comodule over $\chi_{MHTS}$. This is given by the graded vector
space $\Q\langle\langle X_s\rangle\rangle_{s\in S}$ where degree of
$X_s$ is 1 with the coaction:
$$X_{s_1}\cdots X_{s_n}\mapsto \sum X_{t_1}\cdots
X_{t_m}\otimes \I(X_{t_1}\cdots X_{t_m},(X_{s_1}\cdots X_{s_n})').$$
Here $\I(X_{t_1}\cdots X_{t_m},(X_{s_1}\cdots X_{s_n})')$ denotes
framed MHTS given by $\Pi^{uni}({\Bbb A}^1-S;a,b)$ with frames
\begin{align*}
X_{t_1}\cdots X_{t_m}\in& \Pi_{-2m},
\\
(X_{s_1}\cdots X_{s_n})'\in& \mbox{Hom}(\Pi_{-2n},\Q).
\end{align*}
\begin{lem}\label{equivalent}
 The framed MHTS given by $\I(X_{t_1}\cdots X_{t_m},(X_{s_1}\cdots X_{s_n})')$ is equivalent to
$$\sum \prod_{j=0}^m \I^{\mathcal H}(s_{i_j};s_{i_j+1},\dots;s_{i_{j+1}}),$$
where the sum is taken over all the indices $0=i_0<i_1<\dots
i_m<i_{m+1}=n+1$ such that $s_{i_j}=t_j$. Also by convention $s_0=a$
and $s_{n+1}=b$.
\end{lem}
{\bf Proof.} The equivalence is given by the morphism of MHTS: $$
\Pi({\Bbb A}^1-S;a,b)(m)\To \otimes_{j=0}^m \Pi({\Bbb
A}^1-S;t_j,t_{j+1}),$$ defined on the deRham space by:

$$X_{b_1}\cdots X_{b_k} \mapsto \sum \bigotimes_{j=0}^m (X_{b_{i_j+1}}\cdots
X_{b_{i_{j+1}-1}}),$$ where the sum is taken over all indices
$0=i_0<i_1<\dots<i_m<i_{m+1}=k+1$ such that $b_{i_j}=t_j$. The fact
that this gives the desired equivalence is explained in \cite{G2}\S
6.\qed

We can therefore write the co-action $V\To V\otimes \chi_{MHTS}$ as
$$X_{s_1}\cdots X_{s_n}\mapsto \sum
X_{s_{i_1}}\cdots X_{s_{i_k}} \otimes\prod_{j=0}^k \I^{\mathcal
H}(s_{i_j};s_{i_j+1},\dots;s_{i_{j+1}}),$$ where the sum is taken
over ${0=i_0<i_1<\dots<i_k<i_{k+1}=n+1}$ for $k=0,,\dots$ and
$s_0:=a$, $s_{n+1}:=b$. Hence it is natural to make the following
definition:
\begin{defn} Analogously to the Hodge description above, we define the combinatorial motivic torsor of paths, $\Pi^{\bf Com}({\Bbb
A}^1-S;a,b)$ as the graded vector space $\Q\langle\langle
X_s\rangle\rangle_{s\in S}$ where degree of $X_s$ is 1 with the
coaction:
$$X_{s_1}\cdots X_{s_n}\mapsto \sum
X_{s_{i_1}}\cdots X_{s_{i_k}} \otimes(-1)^{n-k} \prod_{j=0}^k
T(s_{i_j};s_{i_j+1},\dots;s_{i_{j+1}}),$$ where the sum is taken
over ${0=i_0<i_1<\dots<i_k<i_{k+1}=n+1}$ for $k=0,,\dots$ and
$s_0:=a$, $s_{n+1}:=b$.
\end{defn}
\begin{thm}
The definition given above, defines a combinatorial MTM over $F$.
We have a morphism in the category of combinatorial MTM over $F$:
$$\Pi^{\bf Com}({\Bbb A}^1-S;a,c)\otimes \Pi^{\bf Com}({\Bbb
A}^1-S;c,b)\To \Pi^{\bf Com}({\Bbb A}^1-S;a,b)$$ which is the
analogue of the path composition.
\end{thm}
{\bf Proof.} It is clearly compatible with the counit, i.e. the
augmentation. To show its compatibility with coproduct we write the
coaction as:
$$X_{\bf s}\mapsto \sum_D X_{P_1}\otimes
(-1)^{n-\deg(P_1)}[t(P_2)|\dots|t(P_m)]$$ where the sum is taken
over all admissible decompositions $D=P_1\cup\cdots\cup P_m$ of
$\bf s$ and $X_{P_1}=X_{s_{i_1}}\cdots X_{s_{i_k}}$ if the
vertices of $P_1$ are $s_{i_1},\dots,s_{i_k}$ and we define
$\deg(P_1)=k$, the number of vertices of $P_1$. If we apply the
coaction one more time we get $$\sum_{D,D'} X_{P_{11}}\otimes
(-1)^{n-\deg(P_{11})}[t(P_{12})|\dots|t(P_{1r})]\otimes
[t(P_2)|\dots|t(P_m)]$$ where $D$ is as above and $D'=P_{11}\cup
\cdots \cup P_{1r}$ is an admissible decomposition of $P_1$. Since
the coproduct is defined by
$$\Delta([a_1|\dots|a_r])=\sum_{s=0}^r
[a_1|\dots|a_s]\otimes[a_{s+1}|\dots|a_r]$$ this proves
compatibility of the coaction with the coproduct. The composition
morphism exists because of the path composition in theorem
\ref{properties} was proved at the level of trees.\qed

\begin{rem}
 Existence of a natural morphism $H^0B(\tilde{\T}_F)\To H^0B(\N)$
 which after composing with the Hodge realization sends
 $T(\bf{a})$ to $\I^{\mathcal H}({\bf a})$ implies the existence of
 the motivic fundamental group in the category of MTM of Bloch and
 K\v ri\v z.
 \end{rem}

\section{Miscellaneous Remarks}

In the previous sections for any sequence
$a_0;a_1,\dots,a_n;a_{n+1}$ of elements of $F$, an element
$t(a_0;a_1,\dots,a_n;a_{n+1})\in \T_F^1(n)$ was constructed. The
crucial property was its differential formula given in Proposition
\ref{tree}. From this one can construct an element
$T(a_0;a_1,\dots,a_n;a_{n+1})\in \h^0B(\T_F)$. The problem of
constructing motivic iterated integrals is that there is no
natural morphism from $\T_F\To \N$, whereas we have a morphism:
$\T_{F,F}\To \N.$ Therefore one needs to have a lifting of
$t(a_0;\dots;a_{n+1})$ to $\T_{F,F}$ that has the same
differential formula. When no non-zero term repeats, such a
lifting exists. Namely by given zero labels to all edges except
the leaves with zero vertex label which get 1 as their edge label.
\\
We were unable to construct a canonical lifting for the general
case. in this section we want to give some partial results that we
obtained in this direction.
\\
\begin{lem}
There is a natural lifting of $t(a_0;a_1,a_2;a_3)$ in $\T_{F,F}$.
\end{lem}

{\bf Proof.} A 3-edge 3-valent tree with root decorated by $a$ and
the leaves decorated by $b$ and $c$, has the following lifting to
$\T_{F,F}$: if $(a,b,c)$ is generic then we label the edges as
before. If $b=c$ then by symmetry it is zero. If $a=b$ and $b\ne
c$ then we label the edge with vertex label $c$ by $b$ and the
other two edges as usual. Similarly if $a=c$ and $b\ne c$, we
label the edge with vertex label $b$ by $c$ and the other two
edges as usual.  It is easy to see that it defines an element of
$\T_{F,F}$ and is compatible with the differential. Using this we
can define the desired lifting of $t(a_0;a_1,a_2;a_3)$.\qed

\begin{lem}
Let ${\bf{a}}=(0;a_1,\dots,a_n;1)$ be a sequence with at most two
$a_i$ equal to $1$ and the rest equal to zero. There is a lifting
of $t({\bf a})$ to $\T_{F,\{0,1\}}$ such that it satisfies the
same differential equation as in proposition \ref{tree}.
\end{lem}

{\bf Proof.} Notice that the terms in the differential are
admissible by the assumption on ${\bf a}$. The only type of trees
that are not admissible in $t({\bf a})$ are those with root
labelled by 1 and the two edges coming off the root edge, one is a
leaf with vertex 1 and the other is an internal edge. To make this
admissible we label this internal edge by 1. This way the
differential of the tree would not be changed and therefore the
same differential equation holds for this lifting.\qed
\begin{eg}
 With this admissible lifting we can apply the method of corollary \ref{bump}.
 We therefore have an element of $\h^0B(\N)$. It satisfies the same
 coproduct formula as for the double zeta framed motive constructed in
 \cite{G1}, therefore its difference with the true double zeta
 motive is in $\mbox{Ext}^1(\Q(0),\Q(n))$, where $n$ is the
 weight of the farmed object. In case $n$ is even and we are over a number field this
 extension is zero and we therefore get the true double zeta
 motive. To do this concisely we define the
 cycles $\rho(n,m,k)$. Here
 $n,m,k$ are non-negative integers. If $n+k>0$, they are given by:
 \begin{align*}
 &\rho(n,m,k)=(-1)^{m+n}\mbox{Alt}(1-\frac{1}{x_1},1-\frac{x_1}{x_2},\dots,
 1-\frac{x_{k+n}}{x_{k+n+1}},
 1-x_{k+n+1},1-\frac{x_{k+1}}{x_{k+n+2}},\\&1-\frac{x_{k+n+2}}{x_{k+n+3}},\dots,1-\frac{x_{k+n+m}}{x_{k+n+m+1}},1-x_{k+n+m+1},
 x_{k+n+m+1},x_{k+n+m},\dots,\widehat{x_{k+1}},\dots,x_1).
\end{align*}
and if $n=k=0$:
\begin{align*}
\rho(0,m,0)=(-1)^m\mbox{Alt}&(1-\frac{1}{x_0},1-x_0,\frac{x_0-x_1}{1-x_1},1-\frac{x_1}{x_2},\dots,1-\frac{x_{m-1}}{x_m},1-x_m,x_{m},\dots,x_1).
\end{align*}
Using this notation the cycle associated to $\zeta(n,m)$ for $n+m$
an even integer, is:
$$\sum_{k\ge 0, l\ge 0} (-1)^{m-k-1}\binom{n+m-k-l-2}{n-l-1}\rho(l,n+m-l-k-2,k).$$
To give an element of $\h^0B(\N)$ representing $\zeta(n,m)$ (for
$n+m$ even) we have:
\begin{align*}
\zeta(n.m)=&\sum_{k\ge 0, l\ge 0}
(-1)^{m-k-1}\binom{n+m-k-l-2}{n-l-1}\rho(l,n+m-l-k-2,k)\\ &
+\sum_{l=1}^{n-1}(-1)^{m+l}\binom{m+n-l-2}{n-l-1}[\zeta(l+1)|\zeta(m+n-l-1)],
\end{align*}
where $\zeta(n)$ is the following cycle:
$$(-1)^{n-1}\mbox{Alt}(1-\frac{1}{x_1},1-\frac{x_1}{x_2},\dots,1-\frac{x_{n-2}}{x_{n-1}},1-x_{n-1},x_{n-1},\dots,x_1).$$
\end{eg}

\end{document}